\pdfoutput=1
\RequirePackage{ifpdf}
\ifpdf 
\documentclass[pdftex]{sigma}
\else
\documentclass{sigma}
\fi

\numberwithin{equation}{section}

\newtheorem{Theorem}{Theorem}[section]
\newtheorem*{Theorem*}{Theorem}
\newtheorem{Corollary}[Theorem]{Corollary}
\newtheorem{Lemma}[Theorem]{Lemma}
\newtheorem{Proposition}[Theorem]{Proposition}
 { \theoremstyle{definition}

 }
\newcommand\gf[2]{\genfrac{}{}{0pt}{}{#1}{#2}}
\usepackage{leftidx}
\usepackage{paralist}
\usepackage{mathrsfs}

\begin{document}

\allowdisplaybreaks

\newcommand{\arXivNumber}{1811.10476}

\renewcommand{\PaperNumber}{060}

\FirstPageHeading

\ShortArticleName{Some Differential Equations for the Riemann $\theta$-Function on Jacobians}

\ArticleName{Some Differential Equations\\ for the Riemann $\boldsymbol{\theta}$-Function on Jacobians}

\Author{Robert WILMS}

\AuthorNameForHeading{R.~Wilms}

\Address{Laboratoire de Math\'ematiques Nicolas Oresme, Universit\'e de Caen--Normandie,\\ BP 5186, 14032 Caen Cedex, France}
\Email{\href{mailto:robert.wilms@unicaen.fr}{robert.wilms@unicaen.fr}}
\URLaddress{\url{https://www.robert-wilms.de/}}

\ArticleDates{Received March 08, 2024, in final form June 26, 2024; Published online July 02, 2024}

\Abstract{We prove some differential equations for the Riemann theta function associated to the Jacobian of a Riemann surface. The proof is based on some variants of a formula by Fay for the theta function, which are motivated by their analogues in Arakelov theory of Riemann surfaces.}

\Keywords{$\theta$-functions; Riemann surfaces; Jacobians; differential equations}

\Classification{14H42; 14G40}

\section{Introduction}
The Riemann $\theta$-function in dimension $g\ge 1$ is given by
\[
\theta\left[\gf{\alpha_1}{\alpha_2}\right](\tau,z)=\sum_{m\in \mathbb{Z}^g}\exp\bigl(\pi {\rm i}\ltrans{(m+\alpha_1)}\tau(m+\alpha_1)+2\pi {\rm i}\ltrans{(m+\alpha_1)}(z+\alpha_2))\bigr),
\]
where $\alpha=(\alpha_1,\alpha_2)\in \bigl(\frac{1}{2}\mathbb{Z}\bigr)^g\times\bigl(\frac{1}{2}\mathbb{Z}\bigr)^g$ is called a theta characteristic, $z\in\mathbb{C}^g$ and $\tau$ is a complex symmetric $g\times g$ matrix with positive definite imaginary part. Fay \cite{Fay73} studied intensively the case, where $\tau$ is the period matrix of a Riemann surface. In particular, he obtained his famous trisecant identity, which can be applied to give solutions of certain differential equations occurring in mathematical physics in terms of $\theta$-functions as in \cite[Chapter~IIIb, Section~4]{Mum84}. Building on Fay's studies, we will derive some new differential equations for the Riemann $\theta$-function associated to Riemann surfaces.

To give the precise statement, let $X$ be a compact and connected Riemann surface of genus~${g\ge 1}$ and $\tau$ a period matrix for $X$. Let us shortly write $\theta(z)=\theta[0](\tau,z)$ and \smash{$\theta_j=\frac{\partial\theta}{\partial z_j}$} and \smash{$\theta_{jk}=\frac{\partial^2\theta}{\partial z_j\partial z_k}$} for the partial derivatives of $\theta$. Further, we define as in \cite{Gua99}
\[J(w_1,\dots,w_g)=\det(\theta_j(w_k)),\qquad w_1,\dots,w_g\in \mathbb{C}^g\]
and as in \cite{dJo08}
\[\eta=\det\begin{pmatrix}
	\theta_{jk} & \theta_j\\ \theta_k&0 \end{pmatrix}.\]
After fixing a symplectic basis of homology $H_1(X,\mathbb{Z})$ and representing curves, we can canonically identify $\mathrm{Pic}_{g-1}(X)\cong \mathrm{Pic}_{0}(X)$ and fix canonical representatives in $\mathbb{C}^g$ for divisors of degree $g-1$. See Section \ref{secfay} for details.
\begin{Theorem}\label{mainthm}
	Let $X$ be a compact and connected Riemann surface of genus $g\ge 1$ and $p_1,\dots,p_g$, $q\in X$ arbitrary points on $X$ in general position. We denote the degree $(g-1)$ divisor $D=\sum_{j=1}^g p_j-q$ and the effective degree $(g-1)$ divisors $D_k=\sum_{j=1}^g p_j-p_k$ for $1\le k\le g$. Then the following equations hold:
\begin{alignat*}{3}
& (i)\quad &&\prod_{k=1}^g\eta(D_k)=(-1)^{g}\left(\frac{J(D_1,\dots,D_g)^{}}{\theta(D)^{g-1}}\right)^{2g}\prod_{j\neq k}^g\frac{\theta(D_j+p_k-q)^2}{\theta(D_j+p_k-p_j)},&\\
& (ii) \quad&& \prod_{k=1}^g\eta((g-1)p_k)=(-1)^{g}\left(\frac{J(D_1,\dots,D_g)^{}}{\theta(D)^{g-1}}\right)^{2g}\prod_{j\neq k}^g\frac{\theta(gp_j-q)^2}{\theta(gp_j-p_k)},&\\
& (iii) \quad && \eta(D_g)^{g-1}=\prod_{k=1}^{g-1}\left(\eta((g-1)p_k)\left(\frac{\theta(D_g+p_k-q)}{\theta(gp_k-q)}\right)^{g-1}\right).&
\end{alignat*}
\end{Theorem}
We will prove the theorem in Section \ref{secfay} by comparing two derived variants of a formula by Fay on $\theta$, the Schottky--Klein prime form $E(\cdot,\cdot)$ and the Brill--Noether matrix. The most difficult problem is to connect the determinant of the Brill--Noether matrix to the determinants $J$ and $\eta$. Especially for $\eta$ this involves ambitious combinatorics. The proof of the theorem is motivated by analogous formulas in Arakelov theory on the normed versions $\|\theta\|$, $\|J\|$ and $\|\eta\|$ of $\theta$, $J$ and~$\eta$. We will discuss them in Section \ref{secarakelov}.
\section{Variations of Fay's formula}\label{secfay}
In this section, we prove Theorem \ref{mainthm}.
Let $X$ be a compact and connected Riemann surface of genus $g\ge 1$ and $A_1,\dots,A_g,B_1,\dots,B_g\in H_1(X,\mathbb{Z})$ a basis of homology such that for all~$j$,~$k$ we have $(A_j,A_k)=(B_j,B_k)=0$ and $(A_j,B_k)=\delta_{jk}$ for the intersection pairings. Further, let $v_1,\dots,v_g\in H^0\bigl(X,\Omega_X^1\bigr)$ be the unique basis of one-forms such that \smash{$\int_{A_k}v_j=\delta_{jk}$} and write~\smash{$v=\ltrans{(v_1,\dots,v_g)}$} for the vector of them. The period matrix $\tau$ of $X$ is given by the entries \smash{$\tau_{jk}=\int_{B_k}v_j$}. It is a complex symmetric $g\times g$ matrix with positive definite imaginary part. Note that $\tau$ is not an invariant of $X$ as it also depends on the choice of the basis of homology $A_1,\dots,A_g,B_1,\dots,B_g\in H_1(X,\mathbb{Z})$. Two different bases of homology as above can be transformed in each other by a symplectic matrix $M=\left(\begin{smallmatrix} A&B\\ C&D\end{smallmatrix}\right)\in\mathrm{Sp}(2g,\mathbb{Z})$. This basis transformation acts on the period matrix by $M\cdot \tau=(A\tau+B)(C\tau+D)^{-1}$.

As in the introduction, we associate to a theta characteristic $\alpha=(\alpha_1,\alpha_2)\in \bigl(\frac{1}{2}\mathbb{Z}\bigr)^{g}\times\bigl(\frac{1}{2}\mathbb{Z}\bigr)^{g}$ and to $\tau$ the theta function
\[\theta[\alpha](z)=\sum_{m\in \mathbb{Z}^g}\exp\bigl(\pi {\rm i}\ltrans{(m+\alpha_1)}\tau(m+\alpha_1)+2\pi {\rm i}\ltrans{(m+\alpha_1)}(z+\alpha_2)\bigr)
\]
on $\mathbb{C}^g$. We shortly write $\theta(z)=\theta[0](z)$ and \smash{$\theta_j=\frac{\partial\theta}{\partial z_j}$} and \smash{$\theta_{jk}=\frac{\partial^2\theta}{\partial z_j\partial z_k}$} for the partial derivatives. A theta characteristic $\alpha$ is called even if $\theta[\alpha](z)$ is an even function or equivalently if $4\ltrans{\alpha_1}\alpha_2$ is even. Otherwise, it is called odd.

Let $\mathrm{Jac}(X)=\mathbb{C}^g/(\mathbb{Z}^g+\tau\mathbb{Z}^g)$ be the Jacobian of $X$. Any divisor $D=\sum_{j=1}^n p_j-\sum_{j=1}^n q_j$ of~degree~$0$ on~$X$ defines an element \smash{$\sum_{j=1}^n \int_{q_j}^{p_j} v$} in $\mathrm{Jac}(X)$. We fix curves $\gamma_1,\dots,\gamma_{2g}$ representing the homology classes $A_1,\dots,A_g$, $B_1,\dots,B_g$. Further, we make the convention that the paths of integration are taken in $X$ cut along these curves. If we assume that the support of $D$ does not intersect $\bigcup_{j=1}^{2g}\gamma_j$ or, more generally, that the points $p_1,\dots, p_n$, $q_1,\dots,q_n$ are in general position, then we obtain a well-defined representative of $D$ in $\mathbb{C}^g$, which we also denote by $D$. There is an isomorphism
\begin{align}\label{picjacisom}
	\varphi\colon\ \mathrm{Pic}_{g-1}(X)\xrightarrow{\cong} \mathrm{Jac}(X),
\end{align}
which sends the divisor $\Theta\subseteq \mathrm{Pic}_{g-1}(X)$ of effective line bundles of degree $g-1$ to the zero divisor of $\theta$, see for example \cite[Corollary II.3.6]{Mum83}. This can be made explicit by the vector of Riemann constants
\[
\varphi\Biggl(\sum_{j=1}^d n_j p_j\Biggr)=\Biggl(\sum_{j=1}^d n_j\int_{p_0}^{p_j} v_k-\frac{\tau_{kk}-1}{2}-\sum_{\gf{j=1}{j\neq k}}^{g}\int_{A_j}v_j(x)\int_{p_0}^x v_k\Biggr)_{1\le k\le g},\]
where $\sum_{j=1}^d n_j=g-1$ and $p_0\in X$ is some fixed base point not lying in $\bigcup_{j=1}^{2g}\gamma_j$. By our convention to take paths of integration in $X$ cut along $\gamma_1,\dots,\gamma_{2g}$, the above vector gives for any degree $g-1$ divisor $D$ with support outside of \smash{$\bigcup_{j=1}^{2g}\gamma_j$} a well-defined representative in $\mathbb{C}^g$, which we also denote by $D$.

Note that $\varphi$ sends line bundles $\mathscr{L}$ with $\mathscr{L}\otimes\mathscr{L}=K_X$, where $K_X$ denotes the canonical bundle of $X$, to theta characteristics in $\frac{1}{2}\mathbb{Z}^g+\frac{1}{2}\tau\mathbb{Z}^g\subseteq \mathbb{C}^g$. Furthermore, $\varphi(\mathscr{L})$ is even if and only if $\dim H^0(X,\mathscr{L})$ is even. We call a~theta~characteristic $\alpha$ non-singular if the corresponding line bundle $\mathscr{L}_\alpha=\varphi^{-1}(\alpha)$ satisfies $\dim H^0(X,\mathscr{L}_\alpha)\le 1$.
For an odd and non-singular theta~characteristic \smash{$\alpha\in\bigl(\frac{1}{2}\mathbb{Z}\bigr)^{2g}$} and $x\in X$, we define the half-order differential $h_\alpha\in H^0(X,\mathscr{L}_\alpha)$ by
\[h_{\alpha}(x)^2=\sum_{j=1}^g\frac{\partial\theta[\alpha]}{\partial z_j}(0)v_j(x).\]
Then the Schottky--Klein prime-form is defined by
\[
E(x,y)=\frac{\theta[\alpha]\bigl(\int_x^yv\bigr)}{h_{\alpha}(x)h_{\alpha}(y)}
\]
for any $x,y\in X$. It satisfies $E(x,y)=-E(y,x)$. Moreover, it has a simple zero on the diagonal in $X\times X$ and it is non-zero for $x\neq y$. We fix a generic effective divisor $\mathcal{A}=\sum_{j=1}^g a_j$ with~${a_j\in X}$ and we define for any $x\in X$ in general position
\[
\sigma_{\mathcal{A}}(x)=\frac{\theta(\mathcal{A}-x)}{\prod_{j=1}^gE(a_j,x)}.
\]

Now let $p_1,\dots,p_g,q\in X$ be arbitrary points of $X$ and set $D=\sum_{j=1}^gp_j-q$ as well as~${D_k=\sum_{j=1}^gp_j-p_k}$.
We will often assume without loss of generality and without mentioning it that these points are in general position.
There exists a constant $c(\mathcal{A})\in \mathbb{C}$ independent of~$q$,~$p_1,\dots, p_g$ such that
\begin{align}\label{fayidentity}
	\theta(D)=c(\mathcal{A})\cdot\frac{\det(v_j(p_k))\cdot\sigma_{\mathcal{A}}(q)\cdot \prod_{j=1}^g E(p_j,q)}{\prod_{j<k}^{g}E(p_j,p_k)\cdot\prod_{j=1}^g\sigma_{\mathcal{A}}(p_j)}.
\end{align}
This was shown by Fay in the proof of \cite[Corollary 2.17]{Fay73} if one sets $\mathcal{A}=\mathcal{B}$ in the notation there. Fay also constructed a multiple $\sigma$ of~$\sigma_{\mathcal{A}}$, which is independent of the choice of $\mathcal{A}$, and he proved a version of equation \eqref{fayidentity} for $\sigma_{\mathcal{A}}$ replaced by~$\sigma$ and a constant $c$ independent of $\mathcal{A}$. But for our purposes, the above equation is sufficient, so we will not discuss the more involved definition of~$\sigma$ here.
We will give two alternative equations for the invariant $c(\mathcal{A})$ by differentiating the theta function. For this purpose, we define as in the introduction
\[J(w_1,\dots,w_g)=\det(\theta_j(w_k))\]
for any $w_1,\dots,w_g\in\mathbb{C}^g$.
\begin{Proposition}\label{propJ}
	The following equation holds:
	\[
J(D_1,\dots,D_g)=(-1)^{\binom{g+1}{2}} c(\mathcal{A})\cdot\frac{\theta(D)^{g-1}\prod_{j<k}^g E(p_j,p_k)}{\sigma_{\mathcal{A}}(q)^{g-1}\prod_{j=1}^g E(p_j,q)^{g-1}}.\]
\end{Proposition}
\begin{proof}
	The proof is motivated by Gu\`ardia's proof \cite{Gua99} of the analogous result in Arakelov theory, see also \eqref{equguardia}.
	For any $1\le k\le g$, let $U_k$ be an open neighbourhood of $p_k$ and $t_k\colon U_k\to \mathbb{C}$ a local coordinate such that $v_j=f_{jk}{\rm d}t_k$ for some functions $f_{jk}$ on $U_k$. By the chain rule, we obtain
	\[\lim_{q\to p_k}\frac{\theta(D)}{t_k(q)-t_k(p_k)}=-\sum_{j=1}^g\theta_j(D_k)f_{jk}(p_k).\]
	The derivative of the prime form can be directly computed to be
	\begin{align}\label{derivE}
		\lim_{q\to p_k} \frac{E(p_k,q)}{t_k(q)-t_k(p_k)}=\frac{1}{{\rm d}t_k}.
	\end{align}
	Putting these equations together, we obtain
	\[\lim_{q\to p_k}\frac{\theta(D)}{E(p_k,q)}=-\sum_{j=1}^g\theta_j(D_k)v_j(p_k).\]
	Applying this to Fay's identity \eqref{fayidentity} and taking the product over all $k$ gives
	\[(-1)^g\prod_{k=1}^g\sum_{j=1}^g\theta_j(D_k)v_j(p_k)=\frac{(-1)^{\binom{g}{2}}\cdot c(\mathcal{A})^g\cdot\det(v_j(p_k))^g}{\prod_{j<k}^{g}E(p_j,p_k)^{g-2}\cdot\prod_{j=1}^g\sigma_{\mathcal{A}}(p_j)^{g-1}}.\]
	
	If the function $q\mapsto\theta\bigl(\sum_{j=1}^g x_j-q\bigr)$ is not identically $0$, it has the zero divisor $\sum_{j=1}^g x_j$. Thus, $\theta(D_k+p_l-q)$ has a zero of order $2$ in $q=p_l$ as a function of $q$ if $l\neq k$. Hence, we get~${\sum_{j=1}^g \theta_j(D_k)v_j(p_l)=0}$ whenever $k\neq l$. Thus, we can rewrite the left-hand side by
	\[\prod_{k=1}^g\sum_{j=1}^g\theta_j(D_k)v_j(p_k)=J(D_1,\dots,D_g)\cdot \det v_j(p_k).\]
	Therefore, we conclude
	\[J(D_1,\dots,D_g)=\frac{(-1)^{\binom{g+1}{2}}c(\mathcal{A})^g\det(v_j(p_k))^{g-1}}{\prod_{j<k}^{g}E(p_j,p_k)^{g-2}\cdot\prod_{j=1}^g\sigma_{\mathcal{A}}(p_j)^{g-1}}.\]
	Now, the proposition follows by combining this formula with Fay's identity \eqref{fayidentity}.
\end{proof}

To obtain the second equation for $c(\mathcal{A})$, we denote as in the introduction
\[\eta=\det\begin{pmatrix}\theta_{jk}&\theta_j\\ \theta_k &0 \end{pmatrix}.\]
\begin{Proposition}\label{propeta} For any $1\le l\le g$, it holds
	\[\eta(D_l)=(-1)^{\binom{g+1}{2}}c(\mathcal{A})^2\prod_{\gf{j=1}{j\neq l}}^{g}\frac{\theta(D_l+p_j-q)}{\sigma_{\mathcal{A}}(q)\sigma_{\mathcal{A}}(p_j)E(p_j,q)^g}.\]
\end{Proposition}
\begin{proof}
	The proof is similar to the proof of Proposition \ref{propJ}, but more involved. As one can deduce the case $g=1$ directly from Proposition \ref{propJ}, we assume $g\ge 2$.
	By symmetry, we may also assume $l=g$.
	We will first prove the proposition for the special case $p_1=\dots=p_{g-1}$.
	
	Let $t\colon U\to \mathbb{C}$ be a local coordinate for an open neighbourhood $U$ of $p_1$ such that $v_j=f_{j}{\rm d}t$ for some functions $f_j$ on $U$.
	The Wronskian determinant is locally given by
	\[W(p_1)=\det\left(\frac{1}{(k-1)!}\frac{d^{k-1}f_j}{{\rm d}t^{k-1}}\right)_{1\le j,k\le g}(p_1)\]
	and defines a non-zero global section $\widetilde{v}=W\cdot({\rm d}t)^{\otimes g(g+1)/2}$ of $\Omega_X^{g(g+1)/2}$.
	It can be directly computed that we have
	\[\lim_{p_g\to p_1}\lim_{p_{g-1}\to p_1}\dots\lim_{p_2\to p_1}\frac{\det(v_j(p_k))}{\prod_{j<k}^g E(p_j,p_k)}=\widetilde{v}(p_1).\]
	Applying this to Fay's identity \eqref{fayidentity}, we obtain
	\begin{align}\label{thetawronsk}
		\theta(gp_1-q)=\frac{c(\mathcal{A})\sigma_{\mathcal{A}}(q) \widetilde{v}(p_1)E(p_1,q)^g}{\sigma_{\mathcal{A}}(p_1)^{g}}.
	\end{align}
	Note that $\theta(gp_1-q)$ vanishes of order $g$ at $p_1$ as a function in $q$. By L'H\^opital's rule and \eqref{derivE}, we deduce
	\[F(p_1):=\lim_{q\to p_1} \frac{\theta(gp_1-q)}{E(p_1,q)^g}=\frac{1}{g!}\left.\frac{\partial^g\theta(gp_1-q)}{\partial q^g}\right|_{q=p_1}{\rm d}t^{\otimes g}.\]
	\begin{Lemma}\label{lemma}
		It holds $F(p_1)^{g+1}=(-1)^{\binom{g+1}{2}}\eta((g-1)p_1)\cdot \widetilde{v}(p_1)^{2}$.
	\end{Lemma}
	\begin{proof}
		The idea of the proof is based on \cite[Section 5]{dJo10}. Let us first shorten notations by setting
		\[\theta_{j_1\dots j_n}=\frac{\partial^n\theta}{\partial z_{j_1}\cdots \partial z_{j_n}} ((g-1)p_1),\qquad f_j^{(k)}=\frac{{\rm d}^{k}f_j}{{\rm d}t^k}(p_1).\]
		Note that $\eta((g-1)p_1)\cdot \widetilde{v}(p_1)^{2}$ is given by ${\rm d}t^{\otimes g(g+1)}$ multiplied by the determinant of the following symmetric matrix:
		\[\begin{pmatrix}\sum_{j,k=1}^g\theta_{jk}f_j f_k & \dots & \sum_{j,k=1}^g\theta_{jk}f_j \frac{f_k^{(g-1)}}{(g-1)!} & \sum_{j=1}^g \theta_j f_j\\
			\vdots&\ddots&\vdots&\vdots\\
			\sum_{j,k=1}^g\theta_{jk}\frac{f_j^{(g-1)}}{(g-1)!} f_k & \dots & \sum_{j,k=1}^g\theta_{jk}\frac{f_j^{(g-1)}}{(g-1)!} \frac{f_k^{(g-1)}}{(g-1)!} & \sum_{j=1}^g \theta_j \frac{f_j^{(g-1)}}{(g-1)!}\\
			\sum_{j=1}^g\theta_{j} f_j & \dots & \sum_{j=1}^g\theta_{j}\frac{f_j^{(g-1)}}{(g-1)!} & 0\end{pmatrix}.\]
		
		Let $s$ be any positive integer. For any vectors $a,b\in \mathbb{Z}^s$ we write $a\le b$ if $a_i\le b_i$ for all~${1\le i\le s}$ and $a<b$ if $a\le b$ and $a\neq b$. Further, we denote $|a|=\sum_{i=1}^sa_i$. For two vectors~${m\in \mathbb{Z}^s}$ and $n\in\mathbb{N}_0^s$, we define
		\[h_{m,n}=\left.\frac{\partial^{|n|}\theta((g-1)p_1+m_1(q_1-p_{1})+\dots+m_s(q_s-p_{1}))}{\partial q_1^{n_1}\cdots\partial q_s^{n_s}}\right|_{q_1=\dots=q_s=p_{1}}.\]
		Then we have $g!\cdot F(p_1)=h_{-1,g}{\rm d}t^{\otimes g}$ and $h_{m,n}=0$ for all $n$ if $m\ge 0$ and $|m|\le g-1$, since the involved theta function is constantly zero as a function in $(q_1,\dots,q_s)$ in these cases. By Fa\`a di Bruno's formula, we can explicitly write
		\begin{align}\label{hmn}
			h_{m,n}=\sum_{0\le l\le n} m^l\left(\prod_{i=1}^s\frac{n_i!}{l_i!}\right)\sum_{k\in [1,\dots,g]^{|l|}}\theta_k\sum_{r_j}\prod_{j=1}^{|l|}\frac{f_{k_j}^{(r_{j}-1)}}{r_{j}!},
		\end{align}
		where the last sum runs over all $|l|$-tuples of positive integers \smash{$r_{1},\dots,r_{|l|}$} which sum to $r_{1}+\dots+r_{l_1}=n_1$, $r_{l_1+1}+\dots+r_{l_1+l_2}=n_2$ and so on. Further, $m^l$ should be read as \smash{$\prod_{j=1}^s m_j^{l_j}$}. Considering $h_{m,n}$ as a multi-degree $n$ polynomial in $m$, we know that this polynomial has to be identically zero for $|n|\le g-1$. In particular, for all $|n|\le g-1$ the coefficients of the monomials of degree $1$ and $2$ in $m$ vanish, which implies
		\[\sum_{j=1}^g\theta_jf_j^{(a)}=\sum_{j,k=1}^g\theta_{jk}f_j^{(b)}f_k^{(c)}=0\]
		for non-negative integers $a$, $b$, $c$ with $a\le g-2$ and $b+c\le g-3$ by choosing ${n=a+1}$ respectively~${n=(b+1,c+1)}$.
		In particular, the determinant of the matrix above is just a~product of~${g+1}$ factors.
		
		Next, we consider the case $s=1$ and $n=g$. We recall that the signed Stirling numbers~${s(g,k)=(-1)^{g-k}\left[\gf{g}{k}\right]}$ are defined as the coefficients of $\prod_{k=0}^{g-1}(X-k)=\sum_{k=0}^g s(g,k)X^k$. Since $h_{m,g}$ is a degree $g$ polynomial vanishing for $m=0,\dots,g-1$, its coefficients have to be given by the signed Stirling numbers multiplied by a common non-zero factor. The factor can be obtained by computing the coefficients of $m$ or $m^2$ in equation \eqref{hmn}
		\begin{align*}
			h_{m,g}&=-\frac{1}{(g-1)!}\sum_{j=1}^g\theta_jf_j^{(g-1)}\sum_{k=1}^g (-1)^{k} \left[\gf{g}{k}\right] m^k\\
			&=\frac{1}{2(g-1)!H_{g-1}}\sum_{r=1}^{g-1}\binom{g}{r}\sum_{j,k=1}^g\theta_{jk}f_j^{(r-1)}f_k^{(g-r-1)}\sum_{k=1}^g (-1)^{k} \left[\gf{g}{k}\right] m^k,
		\end{align*}
		for $H_{k}=\sum_{j=1}^{k}\frac{1}{j}$ the $k$-th harmonic number.
		In particular, we get \[\frac{h_{-1,g}}{g!}=-\sum_{j=1}^g\theta_j\frac{f_j^{(g-1)}}{(g-1)!}=\frac{g}{2H_{g-1}}\sum_{r=1}^{g-1}\frac{1}{r(g-r)}\sum_{j,k=1}^g\theta_{jk}\frac{f_j^{(r-1)}}{(r-1)!}\frac{f_k^{(g-r-1)}}{(g-r-1)!}.\]
		
		We would like to show that the last sum does not depend on the choice of $r$. For this purpose, we consider the case $s=2$ and $n=(r,g-r)$ for some $1\le r\le g-1$. The polynomial $h_{m,n}$ is of multi-degree $n$ and vanishes for all $0\le m<n$. Hence, its coefficients have to be given by the products of signed Stirling numbers $s(g,j)\cdot s(g,k)$ multiplied by a common non-zero factor. The factor can be obtained by considering the coefficient of $m_1 m_2$ in equation \eqref{hmn}
		\begin{align*}
			h_{m,(r,g-r)}=\frac{\sum_{j,k=1}^g\theta_{jk}f_j^{(r-1)}f_k^{(g-r-1)}}{(r-1)!(g-r-1)!}\sum_{j,k\ge 1}(-1)^{j+k}\left[\gf{r}{j}\right]\left[\gf{g-r}{k}\right]m_1^j m_2^k.
		\end{align*}
		Note that the top degree coefficient of $h_{m,(r,g-r)}$ is given by
		\[\sum_{k\in[1,\dots,g]^g}\theta_k\prod_{j=1}^g f_{k_j},\]
		which does not depend on the choice of $r$. Thus,
		\[\sum_{j,k=1}^g\theta_{jk}\frac{f_j^{(r-1)}}{(r-1)!}\frac{f_k^{(g-r-1)}}{(g-r-1)!}\]
		does not depend on the choice of $r$, either. Therefore, the determinant of the matrix above is given by
		\[(-1)^{\binom{g-1}{2}+1}\Biggl(\sum_{j=1}^g\theta_j\frac{f_j^{(g-1)}}{(g-1)!}\Biggr)^2\Biggl(\sum_{j,k=1}^g\theta_{jk}\frac{f_j^{(g-2)}}{(g-2)!}f_k\Biggr)^{g-1}.\]
		Since both expressions in the brackets can be identified with $h_{-1,g}/g!$, we obtain
		\[F(p_1)^{g+1}=\left(\frac{h_{-1,g}}{g!}\right)^{g+1}{\rm d}t^{\otimes g(g+1)}=(-1)^{\binom{g+1}{2}}\eta((g-1)p_1)\cdot\widetilde{v}(p_1)^2,\]
		which proves the lemma.
	\end{proof}

	We continue with the proof of the proposition. Equation \eqref{thetawronsk} implies
	\[F(p_1)=\frac{c(\mathcal{A})\widetilde{v}(p_1)}{\sigma_{\mathcal{A}}(p_1)^{g-1}}.\]
	Applying the lemma, we obtain
	\[\eta((g-1)p_1)=(-1)^{\binom{g+1}{2}}\frac{c(\mathcal{A})^{g+1}\widetilde{v}(p_1)^{g-1}}{\sigma_{\mathcal{A}}(p_1)^{(g-1)(g+1)}}.\]
	Combining this with equation \eqref{thetawronsk}, we obtain
	\[\eta((g-1)p_1)=(-1)^{\binom{g+1}{2}}c(\mathcal{A})^2\left(\frac{\theta(gp_1-q)}{E(p_1,q)^{g}\sigma_{\mathcal{A}}(q)\sigma_{\mathcal{A}}(p_1)}\right)^{g-1},\]
	This proves the proposition in the case $p_1=\dots=p_{g-1}$.
	
	Next we prove the general case. We write \smash{$v_j'\in H^0\bigl(X,\Omega_X^{\otimes 2}\bigr)$} for the two-fold holomorphic differential given locally by $\frac{{\rm d}f_j}{{\rm d}t}({\rm d}t\otimes {\rm d}t)$ for some local coordinate $t$ and $f_j$ such that $v_j=f_j{\rm d}t$. Note that $v_j'$ is independent of the choice of the local coordinate $t$ and hence, it defines a global two-fold holomorphic differential.
	Further, we define for any integer $1\le e\le g-1$ the $g\times g$ matrix
	\[\mathscr{D}_{e}(p_1,\dots,p_{g-1})=\det\begin{pmatrix}
		v_1(p_1)&\dots&v_1(p_{g-1})&v'_1(p_{e})\\
		\vdots&\ddots &\vdots&\vdots\\
		v_g(p_1)&\dots&v_g(p_{g-1})&v'_g(p_{e})\\
	\end{pmatrix}.\]
	We get the following limit:
	\[\lim_{p_g\to p_e}\frac{\det(v_j(p_k))}{E(p_e,p_g)}=\mathscr{D}_e(p_1,\dots,p_{g-1}).\]
	Applying this to Fay's identity \eqref{fayidentity}, we obtain
	\begin{align}\label{equ_theta-dg-pe}
		\theta(D_g+p_e-q)=\frac{c(\mathcal{A}) \mathscr{D}_e(p_1,\dots,p_{g-1})\sigma_{\mathcal{A}}(q) E(p_e,q)\prod_{j=1}^{g-1}E(p_j,q)}{\sigma_{\mathcal{A}}(p_e)\prod_{j\neq e}^{g-1}E(p_j,p_e)\prod_{j<k}^{g-1}E(p_j,p_k)\prod_{j=1}^{g-1}\sigma_{\mathcal{A}}(p_j)},
	\end{align}
	and hence for the product over all $e$
	\begin{align}\label{equ_theta-dg-pe2} \prod_{e=1}^{g-1}\theta(D_g+p_e-q)=\frac{c(\mathcal{A})^{g-1}\prod_{e=1}^{g-1}\mathscr{D}_e(p_1,\dots,p_{g-1})\sigma_{\mathcal{A}}(q)^{g-1}\prod_{j=1}^{g-1}E(p_j,q)^g}{(-1)^{\binom{g-1}{2}}\prod_{j<k}^{g-1}E(p_j,p_k)^{g+1}\prod_{j=1}^{g-1}\sigma_{\mathcal{A}}(p_j)^g}.
	\end{align}
	Since $\theta(D_g+p_e-q)$ vanishes of at least second order at $p_e$ as a function in $q$, we can define
	\begin{align}
		T_e(p_1,\dots,p_{g-1}):={}&\lim_{q\to p_e}\frac{\theta(D_g+p_e-q)}{E(p_e,q)^2} \nonumber\\
		={} &\frac{1}{2}\sum_{j,k=1}^g\theta_{jk}(D_g)v_{j}(p_e)v_{k}(p_e)-\frac{1}{2}\sum_{j=1}^g\theta_j(D_g)v'_{j}(p_e),\label{formulaTe}
	\end{align}
	where the second equality follows by L'H\^opital's rule and Fa\`a di Bruno's formula.
	\begin{Lemma}
		It holds
		\[\prod_{e=1}^{g-1}T_e(p_1,\dots,p_{g-1})^{g+1}=(-\eta(D_g))^{g-1}\prod_{e=1}^{g-1}\mathscr{D}_e(p_1,\dots,p_{g-1})^2.\]
	\end{Lemma}
	\begin{proof}
		The proof is very similar to the proof of Lemma \ref{lemma}.
		We write $v_k=f_{j,k}{\rm d}t$ locally around the point $p_j$ and shorten notations by setting
		\[
\theta_{j_1\dots j_n}=\frac{\partial^n\theta}{\partial z_{j_1}\cdots \partial z_{j_n}} (D_g),\qquad f_{j,k}=f_{j,k}(p_j),\qquad f^{(r)}_{j,k}=\frac{{\rm d}^rf_{j,k}}{{\rm d}t^r}(p_j), \qquad f'_{j,k}=f^{(1)}_{j,k}.\]
		We obtain $\eta(D_g)\mathscr{D}_e(p_1,\dots,p_{g-1})^2$ as ${\rm d}t^{\otimes(2g+2)}$ multiplied by the determinant of the symmetric matrix
		\[M=\sum_{j,k=1}^g\begin{pmatrix}\theta_{jk}f_{1,j} f_{1,k} & \dots & \theta_{jk}f_{1,j} f_{g-1,k} & \theta_{jk}f_{1,j} f'_{e,k} & \frac{1}{g}\theta_j f_{1,j}\\
			\vdots&\ddots&\vdots&\vdots&\vdots\\
			\theta_{jk}f_{g-1,j} f_{1,k} & \dots &\theta_{jk}f_{g-1,j}f_{g-1,k}&\theta_{jk}f_{g-1,j}f'_{e,k}& \frac{1}{g}\theta_j f_{g-1,j}\\
			\theta_{jk}f'_{e,j}f_{1,k}&\dots&\theta_{jk}f'_{e,j}f_{g-1,k}&\theta_{jk}f'_{e,j}f'_{e,k}&\frac{1}{g}\theta_j f'_{e,j}\\
			\frac{1}{g}\theta_{j} f_{1,j} & \dots & \frac{1}{g}\theta_{j}f_{g-1,j} &\frac{1}{g}\theta_{j}f'_{e,j}& 0\end{pmatrix}.\]
		For any $m\in\mathbb{Z}^{g-1}$ and $n\in \mathbb{N}_0^{g-1}$, we define
		\[h'_{m,n}=\frac{\partial^{|n|}\theta(D_g+m_1(q_1-p_1)+\dots+m_{g-1}(q_{g-1}-p_{g-1}))}{\partial q_1^{n_1}\cdots\partial q_{g-1}^{n_{g-1}}}\Big|_{q_1=p_1,\dots,q_{g-1}=p_{g-1}}.\]
		We have $h'_{m,n}=0$ for $0\le m\le (1,\dots,1)$ and all $n$, since the involved theta functions are constantly zero as functions in $q_1,\dots,q_{g-1}$. Fa\`a di Bruno's formula gives
		\[h'_{m,n}=\sum_{0\le l\le n} m^l\left(\prod_{i=1}^{g-1}\frac{n_i!}{l_i!}\right)\sum_{k\in [1,\dots,g]^{l}}\theta_k\sum_{r_{i,j}}\prod_{i=1}^{g-1}\prod_{j=1}^{l_i}\frac{f_{i,k_{i,j}}^{(r_{i,j}-1)}(p_i)}{r_{i,j}!},\]
		where $[1,\dots,g]^{l}$ stands for $\prod_{i=1}^{g-1} [1,\dots,g]^{l_i}$ and the last sum runs over all $|l|$-tuples of positive integers $r_{1,1},\dots,r_{g-1,l_{g-1}}$, which sum to $r_{i,1}+\dots+r_{i,l_i}=n_i$ for every $i$. Considering $h'_{m,n}$ as a~multi-degree $n$ polynomial in $m$, it has to be identically zero for $n\le (1,\dots,1)$. In particular, the coefficients of the monomials of degree $1$ and $2$ in $m$
		\[\sum_{j=1}^g\theta_j f_{a,j}=\sum_{j,k=1}^g \theta_{jk}f_{a,j}f_{b,k}=0\]
		vanish for all $a\neq b$. Hence, the determinant of the matrix $M$ is just a product of $g+1$ factors.
		
		If we set $n(e)=(0,\dots,0,2,0\dots,0)$, where the $2$ occurs at the $e$-th position, the degree $2$ polynomial
		\[h'_{m,n(e)}=m_e^2\sum_{j,k=1}^g \theta_{jk} f_{e,j} f_{e,k}+m_e\sum_{j}\theta_j f'_{e,j}\]
		has to vanish for $m_e=0,1$, which implies $\sum_{j,k=1}^g \theta_{jk} f_{e,j} f_{e,k}=-\sum_{j=1}^g\theta_j f'_{e,j}$ for all $e<g$. Applying this to equation \eqref{formulaTe}, we obtain
		\[T_e(p_1,\dots,p_{g-1})=\sum_{j,k=1}^g \theta_{jk}f_{e,j}f_{e,k}{\rm d}t^{\otimes 2}.\]
		On the other hand, we can compute the determinant of the matrix $M$ as
		\[\eta(D_g)\mathscr{D}_e(p_1,\dots,p_{g-1})^2=-\Biggl(\sum_{j,k=1}^g \theta_{jk}f_{e,j}f_{e,k}\Biggr)^2\Biggl(\prod_{i=1}^{g-1} \sum_{j,k=1}^g \theta_{jk}f_{i,j}f_{i,k}\Biggr){\rm d}t^{\otimes(2g+2)}.\]
		Now the lemma follows by comparing the last two equations after taking the products over all~${e<g}$.
	\end{proof}

	We continue the proof of the proposition.
	Applying the definition in equation \eqref{formulaTe} to equation~\eqref{equ_theta-dg-pe}, we get
	\[T_e(p_1,\dots,p_{g-1})=\frac{c(\mathcal{A})\mathscr{D}_e(p_1,\dots,p_{g-1})}{\prod_{j< k}^{g-1}E(p_j,p_k)\prod_{j=1}^{g-1}\sigma_{\mathcal{A}}(p_j)}.\]
	If we multiply over all $1\le e\le g-1$ and apply the lemma, we get \[\eta(D_g)^{g-1}=(-1)^{(g-1)}\frac{c(\mathcal{A})^{(g-1)(g+1)}\prod_{e=1}^{g-1}\mathscr{D}_e(p_1,\dots,p_{g-1})^{g-1}}{\prod_{j<k}^{g-1}E(p_j,p_k)^{(g-1)(g+1)}\prod_{j=1}^{g-1}\sigma_{\mathcal{A}}(p_j)^{(g-1)(g+1)}}.\]
	Hence, we can conclude by combining with equation \eqref{equ_theta-dg-pe2} \[\eta(D_g)^{g-1}=\Biggl((-1)^{\binom{g+1}{2}}c(\mathcal{A})^2\prod_{j=1}^{g-1}\frac{\theta(D_g+p_j-q)}{\sigma_{\mathcal{A}}(q)\sigma_{\mathcal{A}}(p_j)E(p_j,q)^g}\Biggr)^{g-1}.\]
	This gives the proposition up to a $(g-1)$-th root of unity. But by the special case $p_1=\dots=p_{g-1}$ we know that this root of unity must be $1$.
\end{proof}

\begin{Corollary}\label{cortheta}
	We have the following equalities of meromorphic sections:
	\begin{align*}
		\prod_{j\neq k}^g\frac{\theta(gp_j-q)}{\theta(gp_j-p_k)}=(-1)^{g\binom{g}{2}}\frac{\sigma_{\mathcal{A}}(q)^{g(g-1)}\prod_{j=1}^gE(p_j,q)^{g(g-1)}}{\prod_{j=1}^g\sigma_{\mathcal{A}}(p_j)^{g-1}\prod_{j<k}E(p_j,p_k)^{2g}}=\prod_{j\neq k}^g\frac{\theta(D_j+p_k-q)}{\theta(D_j+p_k-p_j)}.
	\end{align*}
\end{Corollary}
\begin{proof}
	The first equality follows from Proposition \ref{propeta} applied to divisors of the form $(g-1)p_j$ by comparing it for different choices of $q$, namely $q=p_k$ and $q=q$, and multiplying over all~${j\neq k}$. The second equality follows similar, but using the proposition for the divisors $D_j$ instead of $(g-1)p_j$. Note that the $(g-1)$-th root of unity, which may occur, does not depend on the points $p_j$ by continuity. Hence, it cancels out in the quotient.
\end{proof}

\begin{proof}[Proof of Theorem \ref{mainthm}]
	The first two formulas in Theorem \ref{mainthm} are now obtained as combinations of the formulas in Propositions \ref{propJ} and \ref{propeta} and Corollary \ref{cortheta}. The third formula results by comparing the formula in Proposition \ref{propeta} for the divisors $D=\sum_{j=1}^g p_j-q$ and $D=gp_k-q$.
\end{proof}

\section{Analogous results in Arakelov theory}\label{secarakelov}
In this section, we will discuss normed variants of the formulas in Theorem \ref{mainthm} and Section \ref{secfay} in Arakelov theory of Riemann surfaces. We continue the notation from the previous section.
In Arakelov theory, one is interested in canonical norms for sections of line bundles. For the sections $\theta$, $J$ and $\eta$, the norms
\begin{gather*}
\|\theta\|(z)=\det(Y)^{1/4}\exp\bigl(-\pi\ltrans{y}Y^{-1}y\bigr)\cdot |\theta|(z),\\
\|J\|(w_{1},\dots,w_{g})=\det(Y)^{(g+2)/4}\exp\Biggl(-\pi\sum_{k=1}^{g}\ltrans{y_{k}}Y^{-1}y_{k}\Biggr)|J(w_{1},\dots,w_{g})|,\\
\|\eta\|(z)=\det(Y)^{(g+5)/4}\exp\bigl(-\pi(g+1)\ltrans{y}Y^{-1}y\bigr)\cdot|\eta|(z),
\end{gather*}
were given by Faltings \cite[p.~401]{Fal84}, Gu\`ardia \cite[Definition 2.1]{Gua99}, respectively de Jong \cite[Section~2]{dJo08}. Here we denote $Y=\operatorname{Im}(\tau)$, $y=\operatorname{Im}(z)$ and $y_k=\operatorname{Im}(w_k)$ for all $k$.
Arakelov \cite{Ara74} has given a norm for the canonical section of the diagonal bundle $\mathcal{O}_{X^2}(\Delta)$. This norm is the Arakelov--Green function $G(\cdot,\cdot)$ defined by
\[\frac{\partial_{p}\overline{\partial}_{p}}{2\pi {\rm i}}\log G(p,q)^2=\mu(p)-\delta_q(p)\qquad\text{and}\qquad \int_{X}\log G(p,q)\mu(p)=0,\]
where \smash{$\mu=\frac{\rm i}{2g}\sum_{j,k=1}^g \bigl(Y^{-1}\bigr)_{jk}v_jv_k$} is the canonical $(1,1)$-form associated to $X$.

Faltings has given in \cite[p.~402]{Fal84} an analogue of Fay's equation \eqref{fayidentity} for these norms to define his $\delta$-invariant $\delta(X)$, which he used to prove an arithmetic Noether formula. Gu\`ardia \cite[Corollary 3.6]{Gua99} has found an alternative description for $\delta(X)$, which is the following analogue of Proposition~\ref{propJ}:
\begin{align}\label{equguardia}
	\|J\|(D_{1},\dots,D_{g})=e^{-\delta(X)/8}\frac{\|\theta\|(D)^{g-1}\prod_{j<k}G(p_{j},p_{k})}{\prod_{j=1}^gG(p_{j},q)^{g-1}}.
\end{align}
De Jong \cite[Theorem 4.4]{dJo08} has given another formula, which can be expressed as the following analogue of Proposition \ref{propeta}:
\begin{align*}
	\|\eta\|(D_l)=e^{-\delta(X)/4}\prod_{j=1}^{g-1}\frac{\|\theta\|(D_l+p_j-q)}{G(p_j,q)^g},
\end{align*}
see also \cite[equation (2.7)]{Wil17} for this expression.
In the same way as in the proof of Theorem \ref{mainthm}, we can combine both formulas to obtain an analogue of the formulas in the theorem.
\begin{Proposition}
	With the notation as above, we have the following formulas:
	\begin{alignat*}{3}
&(i)\quad &&
			\prod_{k=1}^g\|\eta\|(D_k)=\left(\frac{\|J\|(D_1,\dots,D_g)^{}}{\theta(D)^{g-1}}\right)^{2g}\prod_{j\neq k}^g\frac{\|\theta\|(D_j+p_k-q)^2}{\|\theta\|(D_j+p_k-p_j)},&\\
& (ii)\quad &&
			\prod_{k=1}^g\|\eta\|((g-1)p_k)=\left(\frac{\|J\|(D_1,\dots,D_g)^{}}{\|\theta\|(D)^{g-1}}\right)^{2g}\prod_{j\neq k}^g\frac{\|\theta\|(gp_j-q)^2}{\|\theta\|(gp_j-p_k)},&\\
& (iii)\quad &&
			\|\eta\|(D_g)^{g-1}=\prod_{k=1}^{g-1}\left(\|\eta\|((g-1)p_k)\left(\frac{\|\theta\|(D_g+p_j-q)}{\|\theta\|(gp_j-q)}\right)^{g-1}\right).&
		\end{alignat*}
\end{Proposition}
Of course, the proposition can also be obtained by taking the norms in the formulas of Theorem~\ref{mainthm}. But since the formulas for $\delta(X)$ have been known before, it shows how the proof of the theorem was motivated.

\subsection*{Acknowledgements}
I would like to thank the anonymous referees for their helpful comments and suggestions.
Further, I gratefully acknowledges support from SFB/Transregio 45.


\pdfbookmark[1]{References}{ref}
\LastPageEnding

\end{document}